\documentclass[12pt]{amsart}
\usepackage{xypic}
\usepackage{amscd}
\newtheorem{thm}{Theorem}
\newtheorem{cor}[thm]{Corollary}
\newtheorem{prop}[thm]{Proposition}
\newtheorem{lemma}[thm]{Lemma}

\theoremstyle{remark}
\newtheorem{remark}[thm]{Remark}
\newtheorem{example}[thm]{Example}

\theoremstyle{definition}
\newtheorem{definition}[thm]{Definition}

\numberwithin{equation}{subsection}

\newcommand{\ca}{\mathcal{A}}

\newcommand{\cj}{\mathcal{J}}

\newcommand{\co}{\mathcal{O}}

\newcommand{\im}{\operatorname{Im}}

\newcommand{\Hom}{\operatorname{Hom}}

\newcommand{\Isom}{\operatorname{Isom}}
\newcommand{\isomo}{\overset{\sim}{=}}

\newcommand{\id}{{\mathtt{Id}}}

\newcommand{\shHom}{\underline{\operatorname{Hom}}}
\newcommand{\shEnd}{\underline{\operatorname{End}}}

\newcommand{\shAut}{\underline{\operatorname{Aut}}}
\newcommand{\shIsom}{\underline{\operatorname{Isom}}}

\newcommand{\rk}{\operatorname{rk}}

\newcommand{\pr}{{\mathtt{pr}}}

\newcommand{\DR}{\mathtt{DR}}

\newcommand{\MC}{\operatorname{MC}}

\newcommand{\Def}{\operatorname{Def}}

\newcommand{\Mat}{{\mathtt{Mat}}}

\newcommand{\Der}{{\mathtt{Der}}}

\newcommand{\dbar}{\overline{\partial}}

\DeclareMathOperator{\ad}{ad}

\DeclareMathOperator{\cotr}{cotr} 
\begin{document}

\title{Deformations of Azumaya algebras}

\author[P.Bressler]{Paul Bressler}\thanks{A.~Gorokhovsky was partially  supported by NSF grant DMS-0400342. B.~Tsygan was partially supported by NSF grant DMS-0605030 }
\address{Department of Mathematics, University of Arizona}
\email{bressler@math.arizona.edu}

\author[A.Gorokhovsky]{Alexander Gorokhovsky}
\address{Department of Mathematics, University of Colorado}
\email{Alexander.Gorokhovsky@colorado.edu}

\author[R.Nest]{Ryszard Nest}
\address{Department of Mathematics, University of Copehagen}
\email{rnest@math.ku.dk}

\author[B.Tsygan]{Boris Tsygan}
\address{Department of Mathematics, Northwestern University}
\email{tsygan@math.northwestern.edu}

%\date{\today}
\maketitle

\section{Introduction.}

%\vskip 1cm

%\centerline {\bf Introduction}

In this paper we compute the deformation theory of a special class
of algebras, namely of Azumaya algebras on a manifold ($C^{\infty}$
or complex analytic).

Deformation theory of associative algebras was initiated by
Gerstenhaber in \cite{G}. A deformation of an associative algebra
$A$ over an Artinian ring ${\frak {a}}$ is an ${\frak{a}}$-linear
associative algebra structure on $A\otimes {\frak{a}}$ such that,
for the maximal ideal ${\frak{m}}$ of ${\frak{a}}$, $A\otimes
{\frak{m}}$ is an ideal, and the quotient algebra on $A$ is the
original one. Gerstenhaber showed  that the Hochschild cochain
complex of an associative algebra $A$ has a structure of a
differential graded Lie algebra (DGLA), and that deformations of $A$
over an Artinian ring ${\frak {a}}$ are classified by Maurer-Cartan
elements of the DGLA $C^{\bullet}(A,A)[1]\otimes {\frak {m}}$. A
Maurer-Catan element of a DGLA ${\mathcal L}^{\bullet}$ with the
differential $\delta$ is by definition an element $\lambda$ of
${\mathcal L}^1$ satisfying
\begin{equation} \label{eq:MC}
\delta\lambda +\frac{1}{2}[\lambda, \lambda]=0
\end{equation}
Isomorphic deformations correspond to equivalent Maurer-Cartan
elements, and vice versa.

In subsequent works \cite{GM}, \cite{SS}, \cite{Dr} it was shown
that deformation theories of many other objects are governed by
appropriate DGLAs in the same sense as above. Moreover, if two DGLAs
are quasi-isomorphic, then there is a bijection between the
corresponding sets of equivalence classes of Maurer-Cartan elements.
Therefore, to prove that deformation theories of two associative
algebras are isomorphic, it is enough to construct a chain of
quasi-isomorphisms of DGLAs whose endpoints are the Hochschild
complexes of respective algebras.

This is exactly what is done in this paper. We construct a canonical
isomorphism of deformation theories of two algebras: one is an
Azumaya algebra on a $C^{\infty}$ manifold $X$, the other the
algebra of $ C^{\infty}$ functions on $X$. The systematic study of
the deformation theory of the latter algebra was initiated in
\cite{BFFLS}. In addition to the definition above, it is required
that the multiplication on $A\otimes {\frak{a}}$ be given by
bidifferential expressions. The corresponding Hochschild cochain
complex consists of multi-differential maps $C^{\infty}(X)^{\otimes
n}\to C^{\infty}(X)$. The complete classification of deformations of
$C^{\infty}(X)$ is known from the formality theorem of Kontsevich
\cite{K1}. This theorem asserts that there is a chain of
quasi-isomorphisms connecting the two DGLAs
$C^{\bullet}(C^{\infty}(X), C^{\infty}(X))[1]$ and $\Gamma (X,
\wedge^{\bullet}(TX))[1]$, the latter being the
 DGLA of multivector fields with the Schouten-Nijenhuis bracket.
  For the proofs in the case of a general manifold, cf. also \cite{Do}, \cite{Ha}, \cite{DTT}.

An Azumaya algebra on a manifold $X$ is a sheaf of algebras locally
isomorphic to the algebra of $n\times n$ matrices over the algebra
of functions. Such an algebra determines a second cohomology class
with values in ${\mathcal O}^{\times}_X$ (here, as everywhere in
this paper, in the $C^{\infty}$ case ${\mathcal O}_X$ denotes the
sheaf of smooth functions). This cohomology class $c$
  is necessarily $n$-torsion, i.e. $nc=0$ in $H^2(X,
{\mathcal O}^{\times }_X)$. The main result of this paper is that
there is a chain of quasi-isomorphisms between the Hochschild
complexes of multidifferential cochains of the algebra of functions
and of an arbitrary Azumaya algebra. This chain of
quasi-isomorphisms does depend on some choices but is essentially
canonical. More precisely, what we construct is a canonical
isomorphism in the derived category \cite{H} of the closed model
category \cite{Q} of DGLAs.

Note that in \cite{BGNT} we considered a related problem of
deformation theory. Namely, any cohomology class in $H^2(X,
{\mathcal O}^{\times }_X)$, whether torsion or not, determines an
isomorphism class of a gerbe on $X$, cf. \cite{Br}. A gerbe is a
partial case of an algebroid stack; it is a sheaf of categories on
$X$ satisfying certain properties. In \cite{BGNT} we showed that
deformation theory of algebroid stacks is governed by a certain
DGLA. If an algebroid stack is a gerbe, we constructed a chain of
quasi-isomorphisms between this DGLA and another one, much more
closely related to the Hochschild complex. When the gerbe is an
Azumaya algebra, this latter DGLA is just the Hochschild complex
itself. In ther words, we solved a different deformation problem for
an Azumaya algebra, and got the same answer.

Now let us turn to the case of a complex analytic manifold $X$. In
this case, let ${\mathcal O}_X$ stand for the sheaf of algebras of
holomorphic functions. It is natural to talk about deformations of
this sheaf of algebras. Multidifferential multiholomorphic
Hochschild cochains form a sheaf, and it can be shown that the
corresponding deformation theory is governed by the DGLA of
Dolbeault forms $\Omega ^{0,\bullet}(X,
C^{\bullet>0}({\mathcal{O}}_X, {\mathcal{O}}_X)[1])$. The full DGLA
$\Omega ^{0,\bullet}(X, C^{\bullet}({\mathcal{O}}_X,
{\mathcal{O}}_X)[1])$ governs deformations of ${\mathcal O}_X$ as an
algebroid stack. We show that our construction of a chain of
quasi-isomorphisms can be carried out for this full DGLA in the
holomorphic case, or in a more general case of a real manifold with
a complex integrable distribution. This construction does not seem
to work for the DGLA governing deformations of ${\mathcal O}_X$ as a
sheaf of algebras. As shown in \cite{NT} and \cite{BGNT}, such
deformation theory can be in general more complicated.

Our motivations for studying deformations of gerbes and Azumaya
algebras are the following:

1) A fractional index theorem from \cite{MMS}. The algebra of
pseudodifferential operators which is used there is closely related
to formal deformations of Azumaya algebras.

2) Index theory of Fourier integral operators (FIOs). Guillemin and
Sternberg \cite{GS} have studied FIOs associated to a coisotropic
submanifold of a cotangent bundle. It appears that higher index
theorems for such operators are related to algebraic index theorems
\cite{BNT} for deformations of the trivial gerbe on a symplectic
manifold with an \'{e}tale groupoid. A similar algebraic index
theorem in the holomorphic case should help establish a Riemann-Roch
theorem in the setting of \cite{KS}, \cite{PS}.

3) Dualities between gerbes and noncommutative spaces (\cite{MR1},
\cite{MR2}, \cite{MR3}, \cite{Bl}, \cite{BBP}, \cite{Ka}).

Hochschild and cyclic homology of Azumaya algebras were computed in
\cite{cw94, Schack} and in the more general case of continuous trace
algebras (in the cohomological setting) in \cite{ms06}. Here we
require, however, much more precise statement which involves the
whole Hochschild complex as a differential graded Lie algebra,
rather then just its cohomology groups. It is conjectured in
\cite{Schack}  that algebras which are similar in the sense of
\cite{Schack} have the same deformation theories. Some of our
results can be considered as a verification of this conjecture in
the particular case of Azumaya algebras.

A very recent preprint \cite{Do1} contains results which have a
significant overlap with ours. Namely, a chain of quasi-isomorphisms
similar to ours is established in a partial case when the Azumaya
algebra is an algebra of endomorphisms of a vector bundle. On the
other hand, a broader statement is proven, namely that the above
chain of quasi-isomorphisms extends to Hochschild chain complexes
viewed as DGL modules over DGLAs of Hochschild cochains. This work,
like ours, is motivated by problems of index theory.

\section{Azumaya algebras.}
Let $X$ be a smooth manifold. In what follows we will denote by
$\co_X$ the sheaf of complex-valued $C^{\infty}$ functions on $X$.

\begin{definition}
An \emph{Azumaya algebra} on $X$ is a sheaf of central ${\mathcal
O}_X$-algebras locally isomorphic to $\Mat_n({\mathcal O}_X)$.
\end{definition}

Thus, by definition, the unit map ${\mathcal
O}_X\hookrightarrow{\mathcal A}$ takes values in the center of
${\mathcal A}$.

Let ${\mathcal A}_0 := [{\mathcal A},{\mathcal A}]$ denote the
${\mathcal O}_X$-submodule generated by the image of the commutator
map.

We will now consider ${\mathcal A}$, ${\mathcal A}_0$ and ${\mathcal
O}_X$ as Lie algebras under the commutator bracket. Note that the
bracket on ${\mathcal A}$ takes values in ${\mathcal A}_0$ and the
latter is a Lie ideal in ${\mathcal A}$.

\begin{lemma}
The composition ${\mathcal O}_X\hookrightarrow{\mathcal
A}\to{\mathcal A}/{\mathcal A}_0$ is an isomorphism.
\end{lemma}
\begin{proof}
The issue is local so we may assume that ${\mathcal
A}=\Mat_n({\mathcal O}_X)$ in which case it is well known to be
true.
\end{proof}
\begin{cor}
The map ${\mathcal O}_X\oplus{\mathcal A}_0\to{\mathcal A}$ induced
by the unit map and the inclusion is an isomorphism of Lie algebras.
\end{cor}

\begin{lemma}
The sequence
\[
 0\to{\mathcal O}_X\to{\mathcal A}\stackrel{\ad}{\to}\Der_{{\mathcal
O}_X}({\mathcal A})\to 0
\]
is exact. Moreover, the composition ${\mathcal
A}_0\hookrightarrow{\mathcal A}\stackrel{\ad}{\to}\Der_{{\mathcal
O}_X}({\mathcal A})$ is an isomorphism of Lie algebras.
\end{lemma}

Let ${\mathcal C}({\mathcal A})$ denote the sheaf of (locally
defined) connections on ${\mathcal A}$ with respect to which the
multiplication on ${\mathcal A}$ is horizontal; equivalently, such a
connection $\nabla$ satisfies the Leibniz rule $\nabla(ab)=
\nabla(a)b+a\nabla(b)$ in $\Omega^1_X\otimes_{{\mathcal
O}_X}{\mathcal A}$ for all (locally defined) $a,b\in{\mathcal A}$.
The sheaf ${\mathcal C}({\mathcal A})$ is a torsor under
$\Omega^1_X\otimes_{{\mathcal O}_X}\Der_{{\mathcal O}_X}({\mathcal
A})$.

Any connection $\nabla\in{\mathcal C}({\mathcal A})$ satisfies
$\nabla({\mathcal A}_0)\subset\Omega^1_X\otimes{\mathcal A}_0$.

For $\nabla\in{\mathcal C}({\mathcal A})$ there exists a unique
$\theta = \theta(\nabla)\in\Omega^2_X\otimes{\mathcal A}_0$ such
that $\ad(\theta)=\nabla^2\in\Omega^2_X\otimes\Der_{{\mathcal
O}_X}({\mathcal A})$.

\section{Jets}
Let $\cj_X$ be the sheaf of infinite jets of smooth functions on
$X$. Let $p: \cj_X \to \co_X$ denote the canonical projection.
Suppose now that ${\mathcal A}$ is an Azumaya algebra. Let
$\cj(\ca)$ denote the sheaf of infinite jets of $\ca$. Let $p_{\ca}:
\cj(\ca) \to \ca$ denote the canonical projection. The sheaves
${\mathcal J}({\mathcal A})$ and ${\mathcal A}\otimes_{{\mathcal
O}_X}{\mathcal J}_X$ have canonical structures of sheaves of central
${\mathcal J}_X$-algebras locally isomorphic to $\Mat_n({\mathcal
J}_X)$.

 Let
$\shIsom_0({\mathcal A}\otimes{\mathcal J}_X,{\mathcal J}({\mathcal
A}))$ denote the sheaf of (locally defined) ${\mathcal J}_X$-algebra
isomorphisms ${\mathcal A}\otimes{\mathcal J}_X\to{\mathcal
J}({\mathcal A})$ such that the following diagram is commutative:
\begin{equation*}
\xymatrix{ & \ca \otimes \cj_X \ar[rr]\ar[d]_{\id \otimes p} &&\cj(\ca)\ar[d]^{p_{\ca}}\\
           &\ca \ar[rr]^{\id}&&\ca
          }
\end{equation*}
Similarly denote by  $\shAut_0({\mathcal A}\otimes{\mathcal J}_X)$
  the sheaf of (locally defined) ${\mathcal J}_X$-algebra
automorphisms of ${\mathcal A}\otimes{\mathcal J}_X$ such that the
following diagram is commutative:
\begin{equation*}
\xymatrix{ & \ca \otimes \cj_X \ar[rr]\ar[d]_{\id \otimes p} &&\ca \otimes \cj_X \ar[d]^{\id \otimes p}\\
           &\ca \ar[rr]^{\id}&&\ca
          }
\end{equation*}

\begin{lemma}\label{lemma:isom is a torsor}
The sheaf $\shIsom_0({\mathcal A}\otimes{\mathcal J}_X,{\mathcal
J}({\mathcal A}))$ is a torsor under the sheaf of groups
$\shAut_0({\mathcal A}\otimes{\mathcal J}_X)$.
\end{lemma}
\begin{proof}
Since both ${\mathcal J}({\mathcal A})$ and ${\mathcal
A}\otimes_{{\mathcal O}_X}{\mathcal J}_X$ are locally isomorphic to
$\Mat_n({\mathcal J}_X)$, the sheaf $\shIsom_0({\mathcal
A}\otimes{\mathcal J}_X,{\mathcal J}({\mathcal A}))$ is locally
non-empty, hence a torsor.
\end{proof}

\begin{lemma}
The sheaf of groups $\shAut_0({\mathcal A}\otimes{\mathcal J}_X)$ is
soft.
\end{lemma}
\begin{proof}
Let $\Der_{{\mathcal J}_X, 0}({\mathcal A}\otimes{\mathcal J}_X)$
denote the sheaf of ${\mathcal J}_X$-linear derivations of the
algebra ${\mathcal A}\otimes{\mathcal J}_X$ which reduce to the zero
map modulo ${\mathcal J}_0$. The exponential map
\begin{equation}\label{exponintial map}
\exp : \Der_{{\mathcal J}_X, 0}({\mathcal A}\otimes{\mathcal J}_X)
\to \shAut_0({\mathcal A}\otimes{\mathcal J}_X)\ ,
\end{equation}
$\delta \mapsto \exp(\delta)$, is an isomorphism of sheaves (the
inverse map is given by $a \mapsto \log a = \sum_{n=1}^{\infty}
\frac{(a-1)^n}{n}$). Therefore, it suffices to show that the sheaf
$\Der_{{\mathcal J}_X, 0}({\mathcal A}\otimes{\mathcal J}_X)$ is
soft, but this is clear since it is a module over the sheaf
${\mathcal O}_X$ of $C^\infty$-functions.
\end{proof}
\begin{cor}\label{trivial torsor}
The torsor $\shIsom_0({\mathcal A}\otimes{\mathcal J}_X,{\mathcal
J}({\mathcal A}))$ is trivial, i.e. $\Isom_0({\mathcal
A}\otimes{\mathcal J}_X),{\mathcal J}({\mathcal
A})):=\Gamma(X;\shIsom_0({\mathcal A}\otimes{\mathcal J}_X,{\mathcal
J}({\mathcal A})))\neq\emptyset$.
\end{cor}
\begin{proof}
Since the sheaf of groups $\shAut_0({\mathcal A}\otimes{\mathcal
J}_X)$ is soft we have $H^1(X, \shAut_0({\mathcal A}\otimes{\mathcal
J}_X)=1$ (\cite{DD}, Lemme 22, cf. also \cite{Br}, Proposition
4.1.7). Therefore every $\shAut_0({\mathcal A}\otimes{\mathcal
J}_X)$ torsor is trivial.
\end{proof}

In what follows we will use $\nabla^{can}_{\mathcal E}$ to denote
the canonical flat connection on ${\mathcal J}({\mathcal E})$. A
choice of $\sigma\in\Isom_0({\mathcal A}\otimes{\mathcal
J}_X,{\mathcal J}({\mathcal A}))$ induces the flat connection
$\sigma^{-1}\circ\nabla^{can}_{\mathcal A}\circ\sigma$ on ${\mathcal
A}\otimes{\mathcal J}_X$.

A choice of $\nabla\in{\mathcal C}({\mathcal A})(X)$ give rise to
the connection $\nabla\otimes\id+\id\otimes\nabla^{can}_{\mathcal
O}$ on ${\mathcal A}\otimes{\mathcal J}_X$.

\begin{lemma}\label{lemma:connections compared}
\begin{enumerate}
\item
For any $\sigma\in\Isom_0({\mathcal A}\otimes{\mathcal
J}_X,{\mathcal J}({\mathcal A}))$, $\nabla\in{\mathcal C}({\mathcal
A})(X)$, the difference
\begin{equation}\label{difference of connections}
\sigma^{-1}\circ\nabla^{can}_{\mathcal A}\circ\sigma
-(\nabla\otimes\id + \id\otimes\nabla^{can}_{\mathcal O})
\end{equation}
is  ${\mathcal J}_X$-linear.

\item
There exists a unique $F\in\Gamma(X;\Omega^1_X\otimes{\mathcal
A}_0\otimes{\mathcal J}_X)$ (depending on $\sigma$ and $\nabla$)
such that \eqref{difference of connections} is equal to $\ad(F)$.

\item
Moreover, $F$ satisfies
\begin{equation} \label{formula F}
(\nabla\otimes\id + \id\otimes\nabla^{can}_{\mathcal
O})F+\frac{1}{2}[F, F]+\theta=0
\end{equation}
\end{enumerate}
\end{lemma}
\begin{proof}
We leave the verification of the first claim to the reader. If
follows that \eqref{difference of connections} is a global section
of $\Omega^1_X\otimes\Der_{{\mathcal J}_X}({\mathcal
A}\otimes{\mathcal J}_X)$. Since the map ${\mathcal
A}_0\otimes{\mathcal
J}_X\stackrel{\ad\otimes\id}{\longrightarrow}\Der_{{\mathcal
J}_X}({\mathcal A}\otimes{\mathcal J}_X)$ is an isomorphism the
second claim follows.

 We have
\begin{equation}\label{3.0.4}
\sigma^{-1}\circ\nabla^{can}\circ\sigma = \nabla\otimes\id +
\id\otimes\nabla^{can}_{\mathcal O} + \ad F
\end{equation} where $F\in\Gamma(X;\Omega^1_X\otimes \ca\otimes{\mathcal
J}_{X})$.

Since $(\sigma^{-1}\circ\nabla^{can}\circ\sigma)^2=0$ the element
$(\nabla^{\ca}\otimes\id + \id\otimes\nabla^{can})F+\frac{1}{2}[F,
F]+\theta$  must be central; it also must lie in $\ca_0 \otimes
\cj$. Therefore, it vanishes, i.e. the formula \eqref{formula F}
holds.
\end{proof}

\section{Deformations of Azumaya algebras.}

\subsection{Review of formal deformation theory} Consider a DGLA ${\mathcal L}^{\bullet}$ with the differential $\delta$. A Maurer-Catan element of ${\mathcal L}^{\bullet}$ is by definition an element $\lambda$ of ${\mathcal L}^1$ satisfying
\begin{equation} \label{eq:MC1}
\delta\lambda +\frac{1}{2}[\lambda, \lambda]=0
\end{equation}
Now assume that ${\mathcal L}^{0}$ is nilpotent. Then
${\operatorname{exp}}({\mathcal L}^{0})$ is an algebraic group over
the ring of scalars $k$. This group acts on the set of Maurer-Cartan
elements via
\begin{equation} \label{eq: action on MC}
e^X
(\lambda)={\operatorname{Ad}}_{e^X}(\lambda)-\sum_{n=0}^{\infty}\frac{1}{(n+1)!}{\operatorname{ad}}_X^n(\delta
X)
\end{equation}
Informally,
$$\delta +e^X (\lambda)={\operatorname{Ad}}_{e^X}(\delta +\lambda).$$
We call two Maurer-Cartan elements equivalent if they are in the
same orbit of the above action.
\begin{thm} \label{thm:GM} (\cite{GM})
Let ${\frak{a}}$ be an Artinian algebra with the maximal ideal
${\frak{m}}$. A quasi-isomorphism of DGLAs ${\mathcal L}_1^{\bullet}
\to {\mathcal L}_2^{\bullet}$ induces a bijection between the sets
of equivalence classes of Maurer-Cartan elements of ${\mathcal
L}_1^{\bullet} \otimes {\frak{m}}$ and of ${\mathcal L}_2^{\bullet}
\otimes {\frak{m}}$.
\end{thm}
Given a DGLA  ${\mathcal L}^{\bullet}$ and an Artinian algebra
${\frak{a}}$ with the maximal ideal ${\frak{m}}$, denote by
${\operatorname {MC}}({\mathcal L}^{\bullet})$ the set of
equivalence classes of Maurer-Cartan elements of ${\mathcal
L}^{\bullet} \otimes {\frak{m}}$. \vskip 1cm

\subsection{}
Let $\ca$ be  an Azymaya algebra on $X$. Let $\Def({\mathcal A})$
denote the formal deformation theory of ${\mathcal A}$ as a sheaf of
associative ${\mathbb C}$-algebras, i.e. the groupoid-valued functor
of (commutative) Artin ${\mathbb C}$-algebras which associates to an
Artin algebra ${\mathfrak a}$ the groupoid whose objects are pairs
$(\widetilde{\mathcal A},\phi)$ consisting of a flat ${\mathcal
O}_X\otimes_{\mathbb C}{\mathfrak a}$-algebra $\widetilde{\mathcal
A}$ and a map $\phi : \widetilde{\mathcal A}\to{\mathcal A}$ which
induces an isomorphism $\widetilde{\mathcal A}\otimes_{\mathfrak
a}{\mathbb C}\to{\mathcal A}$. The morphism in $\Def({\mathcal
A})({\mathfrak a})$ are morphisms of ${\mathcal O}_X\otimes_{\mathbb
C}{\mathfrak a}$-algebras which commute with the respective
structure maps to ${\mathcal A}$.

Theorem \ref{main thm:equivalence of Def} is the main result of this
note.

\begin{thm}\label{main thm:equivalence of Def}
Suppose that ${\mathcal A}$ is an Azumaya algebra on $X$. There
exists a canonical equivalence $\Def({\mathcal
A})\isomo\Def({\mathcal O}_X)$.
\end{thm}

The proof of Theorem \ref{main thm:equivalence of Def} will be given
in \ref{proof of main thm}. The main technical ingredient in the
proof is Theorem \ref{main technical ingredient}.

\subsection{Hochschild cochains}\label{ssec:cochains}
Recall that the sheaf $C^n({\mathcal A})$ of Hochschild cochains of
degree $n$ is defined by
\[
C^n({\mathcal A}) := \shHom_{\mathbb C}({\mathcal
A}^{\otimes_{\mathbb C} n},{\mathcal A}) \ .
\]
In the case of Azumaya algebras, we will always consider the subcomplex of local sections of $C^n({\mathcal A})$, i.e. of multidifferential
operators. The link from the deformation theory to the Hochschild theory is provided by the following

\begin{prop}\label{Def equivalent to MC} \cite{G}
For every algebra A, there
exists a canonical equivalence $\Def(A)\isomo\MC(C^\bullet({ A})[1]))$.
\end{prop}
In particular, for an Azumaya algebra ${\mathcal {A}}$ there exists a canonical equivalence $\Def({\mathcal
A})\isomo\MC(\Gamma(X;C^\bullet({\mathcal A})[1]))$.
\begin{proof}
By definition, an element of degree one in the DGLA $C^\bullet({ A})[1]\otimes {\frak m}$ is a map $\lambda: A\otimes A \to A\otimes {\frak m}$. Put $a*b=ab+\lambda (a,b)$. Extend $*$ a binary ${\frak a}$-linear operation on $A\otimes {\frak a}$. Modulo ${\frak m}$, this operation is the multiplication in $A$. Its associativity is equivalent to the Maurer-Cartan equation \eqref{eq:MC}. Two Maurer-Cartan elements are equivalent if and only if there is a map $X:A\to A\otimes {\frak m}$ such that, if one extends it to an ${\frak a}$-linear map $X:A\otimes {\frak m}\to A\otimes {\frak a}$, its exponential ${\operatorname{exp}}{\operatorname{ad}}(X)$ is an isomorphism of the corresponding two associative algebra structures on $A\otimes {\frak a}$. On the other hand, any such isomorphism of algebra structures (which is identical modulo ${\frak m}$ is of the form ${\operatorname{exp}}{\operatorname{ad}}(X)$ for some $X$.
\end{proof}

 Let
\[
C^n({\mathcal J}({\mathcal A})) := \shHom^{cont}_{{\mathcal
O}_X}({\mathcal J}({\mathcal A})^{\otimes_{{\mathcal O}_X}
n},{\mathcal J}({\mathcal A}))
\]
There exists a canonical map $j^\infty : C^n({\mathcal A})\to
C^n({\mathcal J}({\mathcal A}))$ which to a multidifferential
operator associates its linearization. The canonical flat connection
on ${\mathcal J}({\mathcal A})$ induces a flat connection, denoted
$\nabla^{can}$, on $C^n({\mathcal J}({\mathcal A}))$. The de Rham
complex $\DR(C^n({\mathcal J}({\mathcal A})))$ satisfies
$H^i(\DR(C^n({\mathcal J}({\mathcal A})))) = 0$ for $i\neq 0$ while
the map $j^\infty$ induces an isomorphism $C^n({\mathcal A})\isomo
H^0(\DR(C^n({\mathcal J}({\mathcal A}))))$.

The Hochschild differential, denoted $\delta$ and the Gerstenhaber
bracket endow $C^\bullet({\mathcal A})[1]$ (respectively,
$C^\bullet({\mathcal J}({\mathcal A}))[1]$) with a structure of a
DGLA. The connection $\nabla^{can}$ acts by derivations of the
Gerstenhaber bracket on $C^\bullet({\mathcal J}({\mathcal A}))[1]$.
Since it acts by derivations on ${\mathcal J}({\mathcal A})$ the
induced connection on $C^\bullet({\mathcal J}({\mathcal A}))[1]$
commutes with the Hochschild differential. Hence, the graded Lie
algebra $\Omega^\bullet_X\otimes C^\bullet({\mathcal J}({\mathcal
A}))[1]$ equipped with the differential $\nabla^{can}+\delta$ is
DGLA.

\begin{prop}\label{the jet map on cochains is a quism}
The map
\begin{equation}\label{the jet map on cochains}
j^\infty : C^\bullet({\mathcal A})[1] \to\Omega^\bullet\otimes
C^\bullet({\mathcal J}({\mathcal A}))[1]
\end{equation}
is a quasi-isomorphism of DGLA.
\end{prop}
\begin{proof}
It is clear that the map \eqref{the jet map on cochains} is a
morphism of DGLA.

Let $F_i C^\bullet(?) = C^{\geq -i}(?)$. Then, $F_\bullet
C^\bullet(?)$ is a filtered complex and the differential induced on
$Gr^F_\bullet C^\bullet(?)$ is trivial. Consider $\Omega^\bullet_X$
as equipped with the trivial filtration. Then, the map \eqref{the
jet map on cochains} is a morphism of filtered complexes with
respect to the induced filtrations on the source and the target. The
induced map of the associated graded objects a quasi-isomorphism,
hence, so is \eqref{the jet map on cochains}.
\end{proof}

\begin{cor}\label{jets equiv MC}
The map
\[
\MC(\Gamma(X;C^\bullet({\mathcal
A})[1]))\to\MC(\Gamma(X;\Omega^\bullet_X\otimes C^\bullet({\mathcal
J}({\mathcal A}))[1]))
\]
induced by \eqref{the jet map on cochains} is an equivalence.
\end{cor}

\subsection{The cotrace map}
Let $\overline C^n({\mathcal J}_X)$  denote the sheaf of normalized
Hochschild cochains. It is a subsheaf of $C^n({\mathcal J}_X)$ whose
stalks are given by
\begin{multline*}
\overline C^n({\mathcal J}_X)_x=\Hom_{{\mathcal O}_{X,x}}(({\mathcal
J}_{X,x}/{\mathcal O}_{X,x}\cdot 1)^{\otimes n},{\mathcal
J}_{X,x})\subset \\ \Hom{{\mathcal O}_{X,x}}({\mathcal
J}_{X,x}^{\otimes n},{\mathcal J}_{X,x})=C^n({\mathcal J}_X)_x
\end{multline*}
The sheaf $\overline C^\bullet({\mathcal A})[1]$ (respectively,
$\overline C^\bullet({\mathcal J}_X)[1]$ is  actually a sub-DGLA of
$C^\bullet({\mathcal A})[1]$ (respectively, $C^\bullet({\mathcal
J}_X)[1]$ ) and the inclusion map is a quasi-isomorphism.

The flat connection $\nabla^{can}$ preserves $\overline
C^\bullet({\mathcal J}_X)[1]$, and the (restriction to $\overline
C^\bullet({\mathcal O}_X)[1]$ of the) map $j^\infty$ is a
quasi-isomorphism of DGLA.

Consider now the map
\begin{equation}\label{cotrace map}
\cotr:  \overline C^\bullet({\mathcal J}_X)[1]\to
  C^\bullet({\mathcal A}\otimes{\mathcal
J}_X)[1]
\end{equation}
defined as follows:
\begin{equation}\label{define cotrace}
\cotr(D)(a_1\otimes j_1, \dots, a_n \otimes j_n) = a_0\ldots a_n
D(j_1, \ldots, j_n).
\end{equation}
\begin{prop}\label{cotrace is quism}
The map $\cotr$ is a quasiisomorphism of DGLAs.
\end{prop}
\begin{proof}
It is easy to see that $\cotr$ is a morphism of DGLAs. Since the
fact that this is a quasiisomorphism is local it is enough to verify
it when $\ca=\Mat_n({\mathcal O}_X)$. In this case it is a
well-known fact (cf. \cite{Loday}, section 1.5.6).
\end{proof}

\subsection{Comparison of deformation complexes}
Let $\sigma\in\Isom_0({\mathcal A}\otimes{\mathcal J}_X,{\mathcal
J}({\mathcal A}))$, $\nabla\in{\mathcal C}({\mathcal A})$. The
isomorphisms of algebras $\sigma$ induces the isomorphism of DGLA
\[
\sigma_*: C^\bullet({\mathcal A}\otimes{\mathcal J}_X)[1]\to
C^\bullet({\mathcal J}({\mathcal A}))[1]
\]
which is horizontal with respect to the flat connection (induced by)
$\nabla^{can}$ on $C^\bullet({\mathcal J}({\mathcal A}))[1]$ and the
induced flat connection given by \eqref{3.0.4}. Therefore, it
induces the isomorphism of DGLA (the respective de Rham complexes)
\begin{equation}
\sigma_*:   \Omega^\bullet_X\otimes C^\bullet({\mathcal
A}\otimes{\mathcal J}_X)[1]\to \Omega^\bullet_X\otimes
C^\bullet(\cj(\ca) )[1].
\end{equation}
Here, the differential in $ \Omega^\bullet_X\otimes
C^\bullet({\mathcal A}\otimes{\mathcal J}_X)[1]$ is $
\nabla\otimes\id + \id\otimes\nabla^{can} + \ad F +\delta$ and the
differential in $\Omega^\bullet_X\otimes C^\bullet(\cj(\ca) )[1]$ is
$\nabla^{can}+\delta$.

Let $\iota_G$ denote the adjoint action of
$G\in\Gamma(X;\Omega^k_X\otimes C^0({\mathcal A}\otimes{\mathcal
J}_X))$ (recall that $C^0({\mathcal A}\otimes{\mathcal
J}_X)={\mathcal A}\otimes{\mathcal J}_X$). Thus, $\iota_G$ is a map
$\Omega^p_X\otimes C^q({\mathcal A}\otimes{\mathcal
J}_X)\to\Omega^{p+k}_X\otimes C^{q-1}({\mathcal A}\otimes{\mathcal
J}_X)$.

\begin{lemma}\label{commutators}
For any $H, G \in  \Gamma(X;\Omega^{\bullet}_X\otimes C^0({\mathcal
A}\otimes{\mathcal J}_X))$ we have:
\begin{align}
[\delta, \iota_H]&=\ad H\\
[ \ad H, \iota_G]&=\iota_{[H, G]}\\
[\iota_H, \iota_G]&=0\\
[ \nabla\otimes\id + \id\otimes\nabla^{can}_{\co} , \iota_H ]&=
\iota_{(\nabla\otimes\id + \id\otimes\nabla^{can}_{\co})H}
\end{align}
\end{lemma}
\begin{proof}
Direct calculation.
\end{proof}

 Let $\exp(t\iota_F) =
\sum_{k=0}^\infty\frac{1}{n!}t^k(\iota_F)^{\circ k} = \id + t\iota_F
+ \frac{t^2}{2}\iota_F\circ\iota_F + \cdots$. Note that this is a
polynomial in $t$ since $\Omega^p_X=0$ for $p>\dim X$. Since
$\iota_F$ is a derivation, the operation $\exp(\iota_F)$ is an
automorphism of the graded Lie algebra $\Omega^\bullet_X\otimes
C^\bullet({\mathcal A}\otimes{\mathcal J}_X)[1]$. The automorphism
$\exp(\iota_F)$ does not commute with the differential.

\begin{lemma}\label{Adiota}
\begin{multline*}
\exp(\iota_F)\circ (\nabla\otimes\id +
\id\otimes\nabla^{can}_{\mathcal O} + \delta + \ad F)
\circ\exp(-\iota_F) = \\
\nabla\otimes\id + \id\otimes\nabla^{can}_{\mathcal O} + \delta +
\iota_{\theta}.
\end{multline*}
\end{lemma}
\begin{proof}
Consider the following polynomial in $t$:
$p(t)=\exp(t\iota_F)\circ\delta \circ \exp(-t\iota_F)$. Then using
the identities from the Lemma \ref{commutators} we obtain
$p'(t)=\exp(t\iota_F)\circ[\iota_F,\delta] \circ
\exp(-t\iota_F)=-\exp(t\iota_F)\circ \ad F \circ \exp(-t\iota_F)$,
$p''(t)=-\exp(t\iota_F)\circ[\iota_F,\ad F] \circ
\exp(-t\iota_F)=\exp(t\iota_F)\circ \iota_{[F, F]} \circ
\exp(-t\iota_F)$, and $p^{(n)}=0$ for $n \ge 3$. Therefore $p(t)
=\delta -t \ad F +\frac{t^2}{2} \iota_{[F, F]}$. Setting $t=1$ we
obtain
\[
\exp(\iota_F)\circ\delta \circ \exp(-\iota_F)=\delta-\ad F
+\frac{1}{2}\iota_{[F, F]}.
\]
Similarly we obtain
\[
\exp(\iota_F)\circ\ad F \circ \exp(-\iota_F)= \ad F -\iota_{[F, F]}
\]
and
\begin{multline*}
\exp(\iota_F)\circ (\nabla\otimes\id +
\id\otimes\nabla^{can}_{\mathcal O} ) \circ \exp(-\iota_F)=\\
(\nabla\otimes\id + \id\otimes\nabla^{can}_{\mathcal O}) -
\iota_{(\nabla\otimes\id + \id\otimes\nabla^{can}_{\co})F}
\end{multline*}
Adding these formulas up and using the identity \eqref{formula F} we
obtain the desired result.
\end{proof}

\begin{lemma}
The map
\begin{equation}\label{cotrace with forms}
\id\otimes\cotr: \Omega^\bullet_X\otimes\overline
C^\bullet({\mathcal J}_X)[1]\to \Omega^\bullet_X\otimes
C^\bullet({\mathcal A}\otimes{\mathcal J}_X)[1] \ .
\end{equation}
is a quasiisomorphism of DGLA, where the source (respectively, the
target) is equipped with the differential $\nabla^{can}_{\mathcal O}
+ \delta$ (respectively, $\nabla\otimes\id + \id\otimes\nabla^{can}
+ \delta + \iota_{\theta}$).
\end{lemma}
\begin{proof}
It is easy to see that $\id \otimes \cotr$ is a morphism of graded
Lie algebras, which satisfies  $(\nabla \otimes \id +
\id\otimes\nabla^{can}_{\co}) \circ (\id \otimes \cotr)= (\id
\otimes \cotr) \circ \nabla^{can}_{\co} $ and $\delta \circ (\id
\otimes \cotr)= (\id \otimes \cotr) \circ \delta $. Since the domain
of $(\id \otimes \cotr)$ is the normalized complex, we also have
$\iota_{\theta} \circ(\id \otimes \cotr)=0$. This implies that $(\id
\otimes \cotr)$ is a morphism of DGLA.

Introduce filtration on $\Omega^{\bullet}_X$ by
$F_i\Omega^\bullet_X=\Omega^{\geq -i}_X$ (the ``stupid" filtration)
and consider the complexes $\overline C^\bullet({\mathcal J}_X)[1]$
and $C^\bullet({\mathcal A}\otimes{\mathcal J}_X)[1]$ equipped with
the trivial filtration. The map \eqref{cotrace with forms} is a
morphism of filtered complexes with respect to the induced
filtrations on the source and the target. The differentials induced
on the associated graded complexes are $\delta$ (or, more precisely,
$\id\otimes\delta$) and the induced map of the associated graded
objects is $\id\otimes\cotr$ which is a quasi-isomorphism in virtue
of Proposition \ref{cotrace is quism}. Therefore, the map
\eqref{cotrace with forms} is a quasiisomorphism as claimed.
\end{proof}

\begin{prop}\label{prop:quism of DGLA}
For a any choice of $\sigma\in\Isom_0({\mathcal A}\otimes{\mathcal
J}_X,{\mathcal J}({\mathcal A}))$, $\nabla\in{\mathcal C}({\mathcal
A})$, the composition $\Phi_{\sigma,\nabla}:=
\sigma_*\circ\exp(\iota_{F_{\sigma,\nabla}})\circ(\id\otimes\cotr)$
(where $F$ is as in Lemma \ref{lemma:connections compared}),
\begin{equation}\label{map:the comparison map}
\Phi_{\sigma,\nabla}: \Omega^\bullet_X\otimes\overline
C^\bullet({\mathcal J}_X)[1] \to\Omega^\bullet_X\otimes
C^\bullet({\mathcal J}({\mathcal A}))[1]
\end{equation}
is a quasi-isomorphism of DGLA.
\end{prop}
\begin{proof}
This is a direct consequence of the Lemmata \ref{Adiota},
\ref{cotrace with forms}.
\end{proof}

\subsection{Independence of choices}
According to Proposition \ref{prop:quism of DGLA}, for any choice of
$\sigma\in\Isom_0({\mathcal A}\otimes{\mathcal J}_X,{\mathcal
J}({\mathcal A}))$, $\nabla\in{\mathcal C}({\mathcal A})$ we have a
quasi-isomorphism of DGLA $\Phi_{\sigma,\nabla}$.
\begin{prop}\label{prop:independence of choices}
The image of $\Phi_{\sigma,\nabla}$ in the derived category is
independent of the choices made.
\end{prop}
\begin{proof}
For $i=0, 1$ suppose given $\sigma_i\in\Isom_0({\mathcal
A}\otimes{\mathcal J}_X,{\mathcal J}({\mathcal A}))$,
$\nabla_i\in{\mathcal C}({\mathcal A})$. Let $\Phi_i =
\Phi_{\sigma_i,\nabla_i}$. The goal is to show that $\Phi_0=\Phi_1$
in the derived category.

There is a unique $\theta_i\in\Gamma(X;\Omega^2_X\otimes{\mathcal
A}_0\otimes{\mathcal J}_X))$ such that $\nabla_i^2=\ad(\theta_i)$.
By Lemma \ref{lemma:connections compared} there exist unique
$F_i\in\Gamma(X;\Omega^1_X\otimes{\mathcal A}_0\otimes{\mathcal
J}_X)$ such that
\[
\ad(F_i)=\sigma_i^{-1}\circ\nabla^{can}_{\mathcal A}\circ\sigma_i
-(\nabla_i\otimes\id + \id\otimes\nabla^{can}_{\mathcal O})
\]

Let $I:=[0,1]$. Let $\epsilon_i : X\to I\times X$ denote the map
$x\mapsto (i,x)$. Let $\pr:I\times X\to X$ denote the projection on
the second factor.

It follows from Lemma \ref{lemma:isom is a torsor}, the isomorphism
\eqref{exponintial map} and the isomorphism ${\mathcal
A}_0\otimes{\mathcal J}_X\to\Der_{{\mathcal J}_{X,0}}({\mathcal
A}\otimes{\mathcal J}_X)$ that there exists a unique
$f\in\Gamma(X;{\mathcal A}_0\otimes{\mathcal J}_X)$ such that
$\sigma_1=\sigma_0\circ\exp(\ad(f))$. For $t\in I$ let
$\sigma_t=\sigma_0\circ\exp(\ad(tf))$. Let $\widetilde\sigma$ denote
the isomorphism $\pr^*({\mathcal A}\otimes{\mathcal
J}_X)\to\pr^*{\mathcal J}({\mathcal A})$ which restricts to
$\sigma(t)$ on $\{t\}\times X$. In particular,
$\sigma_i=\epsilon^*_i(\widetilde\sigma)$.

For $t\in I$ let $\nabla_t = t\nabla_0+(1-t)\nabla_1$. Let
$\widetilde\nabla$ denote the connection on $\pr^*{\mathcal A}$
which restricts to $\nabla_t$ on $\{t\}\times X$. In particular,
$\nabla_i=\epsilon^*_i(\widetilde\nabla)$.

Suppose that
\begin{enumerate}
\item $\widetilde\sigma :
\pr^*({\mathcal A}\otimes{\mathcal J}_X)\to\pr^*{\mathcal
J}({\mathcal A})$ is an isomorphism of $\pr^*{\mathcal
J}_X$-algebras which reduces to the identity map on $\pr^*{\mathcal
A}$ modulo ${\mathcal J}_0$ and satisfies
$\sigma_i=\epsilon^*_i(\widetilde\sigma)$;

\item $\widetilde\nabla\in{\mathcal C}(\pr^*{\mathcal A})$ satisfies
$\nabla_i=\epsilon^*_i(\widetilde\nabla)$.
\end{enumerate}
(Examples of such are constructed above.)

Then, there exists a unique $\widetilde{F}\in\Gamma(I\times
X;\Omega^1_{I\times X}\otimes\pr^*({\mathcal A}_0\otimes{\mathcal
J}_X))$ such that
\[
\widetilde\sigma^{-1}\circ\pr^*(\nabla^{can}_{\mathcal
A})\circ\widetilde\sigma = \widetilde\nabla\otimes\id +
\id\otimes\pr^*(\nabla^{can}_{\mathcal O}) + \ad(\widetilde{F})
\]
It follows from the uniqueness that $\epsilon^*(\widetilde{F})=F_i$.

There exists a unique $\widetilde\theta\in\Gamma(I\times
X;\Omega^2_{I\times X}\otimes\pr^*({\mathcal A}_0\otimes{\mathcal
J}_X))$ such that $\widetilde\nabla^2=\ad(\widetilde\theta)$. It
follows from the uniqueness that
$\epsilon_i^*(\widetilde\theta)=\theta_i$.

The composition
$\widetilde\Phi:=\widetilde\sigma_*\circ\exp(\iota_{\widetilde{F}})\circ
(\id\otimes\cotr)$,
\begin{equation}\label{interpolation}
\widetilde\Phi:\Omega^\bullet_{I\times X}\otimes\pr^*\overline
C^\bullet({\mathcal J}_X))[1]\to\Omega^\bullet_{I\times
X}\otimes\pr^*C^\bullet({\mathcal J}({\mathcal A}))[1]
\end{equation}
is a quasi-isomorphism of DGLA, where the differential on source
(respectively, target) is $\pr^*(\nabla^{can}_{\mathcal O})+\delta$
(respecively, $\pr^*(\widetilde\nabla^{can}_{\mathcal A}) + \delta +
\iota_{\widetilde\theta}$).

The map \eqref{interpolation} induces the map of direct images under
the projection $\pr$

\begin{equation}\label{direct image of interpolation}
\widetilde\Phi:\pr_*\Omega^\bullet_{I\times X}\otimes\overline
C^\bullet({\mathcal J}_X)[1]\to\pr_*\Omega^\bullet_{I\times
X}\otimes C^\bullet({\mathcal J}({\mathcal A}))[1]
\end{equation}
which is a quasi-isomorphism (since all higher direct images
vanish).

The pull-back of differential forms $\pr^* :
\Omega^\bullet\to\pr_*\Omega^\bullet_{I\times X}$ is a
quasi-isomorphism of commutative DGA inducing the quasi-isomorphism
of DGLA $\pr^*\otimes\id:\Omega^\bullet_X\otimes\overline
C^\bullet({\mathcal J}_X)[1]\to\pr_*\Omega^\bullet_{I\times
X}\otimes\overline C^\bullet({\mathcal J}_X)[1]$ (with differentials
$\nabla^{can}_{\mathcal O}+\delta$ and $\pr^*(\nabla^{can}_{\mathcal
O})+\delta$ respectively).

The diagram of quasi-isomorphisms of DGLA
\begin{equation}
 \xymatrix{ &&&&\Omega^\bullet_X\otimes C^\bullet({\mathcal J}({\mathcal A}))[1]  \\
 &\Omega^\bullet_X\otimes C^\bullet(\cj_X)[1]\ar[rrr]^{\widetilde\Phi\circ(\pr^*\otimes\id)} \ar[rrru]^{\Phi_0}
 \ar[rrrd]_{\Phi_1}&&&\pr_*\Omega^\bullet_{I\times X}\otimes
 C^\bullet({\mathcal J}({\mathcal A}))[1]\ar[u]_{\epsilon_0^*\otimes \id}
 \ar[d]^{\epsilon_1^*\otimes \id}\\
 && & & \Omega^\bullet_X\otimes C^\bullet({\mathcal J}({\mathcal
 A}))[1]
         }
\end{equation}
is commutative. Since $\epsilon^*_i \circ \pr^* =\id$ for $i=0$, $1$
and $\pr^* \otimes \id$ is a quasiisomorphism,  $\epsilon^*_0
\otimes \id$ and $\epsilon_1^* \otimes \id$ represent the same
morphism in the derived category. Hence  so do $\Phi_0$ and
$\Phi_1$.
\end{proof}

\subsection{The main technical ingredient}
To each pair $(\sigma,\nabla)$ with $\sigma\in\Isom_0({\mathcal
A}\otimes{\mathcal J}_X,{\mathcal J}({\mathcal A}))$ and
$\nabla\in{\mathcal C}({\mathcal A})$ we associated the
quasi-isomorphism of DGLA \eqref{map:the comparison map}
(Proposition \ref{prop:quism of DGLA}). According to Proposition
\ref{prop:independence of choices} all of these give rise to the
same isomorphism in the derived category. We summarize these
findings in the following theorem.

\begin{thm}\label{main technical ingredient}
Suppose that ${\mathcal A}$ is an Azumaya algebra on $X$. There
exists a canonical isomorphism in the derived category of DGLA
$\Omega^\bullet_X\otimes\overline C^\bullet({\mathcal
J}_X)[1]\isomo\Omega^\bullet_X\otimes C^\bullet({\mathcal
J}({\mathcal A}))[1]$.
\end{thm}

\subsection{The proof of Theorem \ref{main thm:equivalence of
Def}}\label{proof of main thm} The requisite equivalence is the
composition of the equivalences
\begin{eqnarray*}
\Def({\mathcal O}_X) & \isomo & \MC(\Gamma(X;\overline
C^\bullet({\mathcal
O}_X)[1])) \\
& \isomo & \MC(\Gamma(X;\Omega^\bullet_X\otimes\overline
C^\bullet({\mathcal
J}_X)[1])) \\
& \isomo & \MC(\Gamma(X;\Omega^\bullet_X\otimes C^\bullet({\mathcal
J}({\mathcal A}))[1])) \\
& \isomo & \MC(\Gamma(X;C^\bullet({\mathcal A})[1])) \\
& \isomo & \Def({\mathcal A})
\end{eqnarray*}
The first and the last equivalences are those of Proposition
\ref{Def equivalent to MC}, the second and the fourth are those of
Corollary \ref{jets equiv MC}, and the third one is induced by the
canonical isomorphism in the derived category of Theorem \ref{main
thm:equivalence of Def}.

\section{Holomorphic case}

\subsection{Complex distributions}
Let ${\mathcal T}_X$ denote the sheaf of \emph{real valued} vector
fields on $X$ and let ${\mathcal T}_X^{\mathbb C}:={\mathcal
T}_X\otimes_{\mathbb R}{\mathbb C}$.

\begin{definition}
A \emph{(complex) distribution} on $X$ is a sub-bundle of ${\mathcal
T}_X^{\mathbb C}$\footnote{A sub-bundle is an ${\mathcal
O}_X$-submodule which is a direct summand locally on $X$}
\end{definition}

For a distribution ${\mathcal P}$ on $X$ we denote by ${\mathcal
P}^\perp\subseteq\Omega^1_X$ the annihilator of ${\mathcal P}$ (with
respect to the canonical duality pairing).

\begin{definition}\label{definition of integrability}
A distribution ${\mathcal P}$ of rank $r$ on $X$ is called
\emph{integrable} if, locally on $X$, there exist functions
$f_1,\ldots , f_r\in{\mathcal O}_X$ such that $df_1,\ldots , df_r$
form a local frame for ${\mathcal P}^\perp$.
\end{definition}

\begin{lemma} An integrable distribution is \emph{involutive}, i.e. it is a
Lie subalgebra of ${\mathcal T}_X^{\mathbb C}$ (with respect to the
Lie bracket of vector fields).
\end{lemma}

\subsection{Differential calculus on $X/{\mathcal P}$}
In this subsection we briefly review relevant definitions and
results of the differential calculus in the presence of integrable
complex distribution. We refer the reader to \cite{Kostant},
\cite{Rawnsley} and \cite{FW} for details and proofs. Suppose that
${\mathcal P}$ is an integrable distribution. Let
$F_\bullet\Omega^\bullet_X$ denote the filtration by the powers of
the differential ideal generated by ${\mathcal P}^\perp$, i.e.
$F_{-i}\Omega^j_X = \bigwedge^i{\mathcal
P}^\perp\wedge\Omega^{j-i}_X\subseteq\Omega^j_X$. Let $\dbar$ denote
the differential in $Gr^F\Omega^\bullet_X$. The wedge product of
differential forms induces a structure of a commutative DGA on
$(Gr^F\Omega^\bullet_X,\dbar)$.

\begin{lemma}
The complex $Gr^F_{-i}\Omega^\bullet_X$ satisfies
$H^j(Gr^F_{-i}\Omega^\bullet_X,\dbar)=0$ for $j\neq i$.
\end{lemma}

Let $\Omega^i_{X/{\mathcal P}}:=
H^i(Gr^F_{-i}\Omega^\bullet_X,\dbar)$, ${\mathcal O}_{X/{\mathcal
P}}:=\Omega^0_{X/{\mathcal P}}$. We have ${\mathcal O}_{X/{\mathcal
P}}:=\Omega^0_{X/{\mathcal P}}\subset{\mathcal O}_X$,
$\Omega^1_{X/{\mathcal P}}\subset{\mathcal
P}^\perp\subset\Omega^1_X$ and, more generally,
$\Omega^i_{X/{\mathcal P}}\subset\bigwedge^i{\mathcal
P}^\perp\subset\Omega^i_X$. The wedge product of differential forms
induces a structure of a graded-commutative algebra  on
$\Omega^\bullet_{X/{\mathcal P}}:=\oplus_i\Omega^i_{X/{\mathcal
P}}[-i]=H^\bullet(Gr^F\Omega^\bullet_X,\dbar)$.

If $f_1,\ldots , f_r$ are as in \ref{definition of integrability},
then $\Omega^1_{X/{\mathcal P}}=\sum_{i=1}^r{\mathcal
O}_{X/{\mathcal P}}\cdot df_i$, in particular,
$\Omega^1_{X/{\mathcal P}}$ is a locally free ${\mathcal
O}_{X/{\mathcal P}}$-module. Moreover, the multiplication induces an
isomorphism $\bigwedge^i_{{\mathcal O}_{X/{\mathcal
P}}}\Omega^1_{X/{\mathcal P}}\to\Omega^i_{X/{\mathcal P}}$.

The de Rham differential $d$ restricts to the map $d :
\Omega^i_{X/{\mathcal P}}\to\Omega^{i+1}_{X/{\mathcal P}}$ and the
complex $\Omega^\bullet_{X/{\mathcal P}}:=(\Omega^i_{X/{\mathcal
P}},d)$ is a commutative DGA. Moreover, the inclusion
$\Omega^\bullet_{X/{\mathcal P}}\to\Omega^\bullet_X$ is a
quasi-isomorphism.
\begin{example}\label{real case}
Suppose that $\overline{\mathcal P}={\mathcal P}$. Then, ${\mathcal
P}={\mathcal D}\otimes_{\mathbb R}{\mathbb C}$, where ${\mathcal D}$
a subbundle of ${\mathcal T}_X$. Then, ${\mathcal D}$ is an
integrable real distribution which defines a foliation on $X$ and
$\Omega^\bullet_{X/{\mathcal P}}$ is the complex of basic forms.
\end{example}

\begin{example} Suppose that $\bar{\mathcal{P}}
\cap \mathcal{P} =0$ and $\bar{\mathcal{P}} \oplus \mathcal{P}
={\mathcal T}_X^{\mathbb C}$. In this case again $\mathcal{P}$ is
integrability is equivalent to involutivity, by Newlander-Nirenberg
theorem. $\Omega^\bullet_{X/{\mathcal P}}$ in this case is a
holomorphic de Rham complex.
\end{example}

\subsection{$\overline\partial$-operators}
Suppose that ${\mathcal E}$ is a vector bundle on $X$, i.e. a
locally free ${\mathcal O}_X$-module of finite rank. A connection
along ${\mathcal P}$ on ${\mathcal E}$ is, by definition, a map
$\nabla^{\mathcal P}: {\mathcal E}\to\Omega^1_X/{\mathcal
P}^\perp\otimes_{{\mathcal O}_X}{\mathcal E}$ which satisfies the
Leibniz rule $\nabla^{\mathcal P}(fe)=f\nabla^{\mathcal
P}(e)+\overline\partial f\cdot e$. A connection along ${\mathcal P}$
gives rise to the ${\mathcal O}_X$-linear map $\nabla^{\mathcal
P}_{(\bullet)}:{\mathcal P}\to\shEnd_{\mathbb C}({\mathcal E})$
defined by ${\mathcal P}\ni\xi\mapsto\nabla^{\mathcal P}_\xi$, with
$\nabla^{\mathcal P}_\xi(e)=\iota_\xi\nabla^{\mathcal P}(e)$.

Conversely, an ${\mathcal O}_X$-linear map $\nabla^{\mathcal
P}_{(\bullet)}:{\mathcal P}\to\shEnd_{\mathbb C}({\mathcal E})$
which satisfies the Leibniz rule $\nabla^{\mathcal P}_\xi(fe)=
f\nabla^{\mathcal P}_\xi(e)+\overline\partial f\cdot e$ determines a
connection along ${\mathcal P}$. In what follows we will not
distinguish between the two avatars of a connection along ${\mathcal
P}$ described above. Note that, as a consequence of the
$\dbar$-Leibniz rule a connection along ${\mathcal P}$ is ${\mathcal
O}_{X/{\mathcal P}}$-linear.

A connection along ${\mathcal P}$ on ${\mathcal E}$ is called flat
if the corresponding map $\nabla^{\mathcal P}_{(\bullet)}:{\mathcal
P}\to\shEnd_{\mathbb C}({\mathcal E})$ is a morphism of Lie
algebras. We will refer to a flat connection along ${\mathcal P}$ on
${\mathcal E}$ as a $\overline\partial$-operator on ${\mathcal E}$.

\begin{example}
The differential $\overline\partial$ in $Gr^F\Omega^\bullet_X$ gives
rise to canonical $\dbar$-operators on $\bigwedge^i{\mathcal
P}^\perp$, $i=0, 1,\ldots$.
\end{example}

\begin{example}\label{hol vector fields}
The adjoint action of ${\mathcal P}$ on ${\mathcal T}_X^{\mathbb C}$
preserves ${\mathcal P}$, hence descends to an action of the Lie
algebra ${\mathcal P}$ on ${\mathcal T}_X^{\mathbb C}/{\mathcal P}$.
The latter action is easily seen to be a connection along ${\mathcal
P}$, i.e. a canonical $\dbar$-operator on ${\mathcal T}_X^{\mathbb
C}/{\mathcal P}$ which is easily seen to coincide with the one
induced on ${\mathcal T}_X^{\mathbb C}/{\mathcal P}$ via the duality
pairing between the latter and ${\mathcal P}^\perp$. In the
situation of Example \ref{real case} this connection is known as the
Bott connection.
\end{example}

\begin{example}\label{example of dbar}
Suppose that ${\mathcal F}$ is a locally free ${\mathcal
O}_{X/{\mathcal P}}$-module of finite rank. Then, ${\mathcal
O}_X\otimes_{{\mathcal O}_{X/{\mathcal P}}}{\mathcal F}$ is a
locally free ${\mathcal O}_X$-module of rank $\rk_{{\mathcal
O}_{X/{\mathcal P}}}{\mathcal F}$ and is endowed in a canonical way
with the $\dbar$-operator, namely, $\id\otimes\dbar$.
\end{example}

A connection on ${\mathcal E}$ along ${\mathcal P}$ extends uniquely
to a derivation of the graded $Gr^F_0\Omega^\bullet_X$-module
$Gr^F_0\Omega^\bullet_X\otimes_{{\mathcal O}_X}{\mathcal E}$. A
$\dbar$-operator $\dbar_{\mathcal E}$ satisfies $\dbar_{\mathcal
E}^2=0$. The complex $(Gr^F_0\Omega^\bullet_X\otimes_{{\mathcal
O}_X}{\mathcal E}, \dbar_{\mathcal E})$ is referred to as the
(corresponding) $\dbar$-complex. Since $\dbar_{\mathcal E}$ is
${\mathcal O}_{X/{\mathcal P}}$-linear, the sheaves $H^i({\mathcal
E},\dbar_{\mathcal E})$ are ${\mathcal O}_{X/{\mathcal P}}$-modules.

\begin{lemma}\label{dbar lemma}
Suppose that ${\mathcal E}$ is a vector bundle and $\dbar_{\mathcal
E}$ is a $\dbar$-operator on ${\mathcal E}$. Then,
$H^i(Gr^F_0\otimes_{{\mathcal O}_X}{\mathcal E},\dbar_{\mathcal
E})=0$ for $i\neq 0$, i.e. the $\dbar$-complex is a resolution of
$\ker(\dbar_{\mathcal E})$. Moreover, $\ker(\dbar_{\mathcal E})$ is
locally free over ${\mathcal O}_{X/{\mathcal P}}$ of rank
$\rk_{{\mathcal O}_X}{\mathcal E}$ and the map ${\mathcal
O}_X\otimes_{{\mathcal O}_{X/{\mathcal P}}}\ker(\dbar_{\mathcal
E})\to{\mathcal E}$ (the ${\mathcal O}_X$-linear extension of the
inclusion $\ker(\dbar_{\mathcal E})\to{\mathcal E}$) is an
isomorphism.
\end{lemma}

\begin{remark} In the notations of Example \ref{example of dbar} and Lemma
\ref{dbar lemma}, the assignments ${\mathcal F}\mapsto({\mathcal
O}_X\otimes_{{\mathcal O}_{X/{\mathcal P}}}{\mathcal
F},\id\otimes\dbar)$ and $({\mathcal E},\dbar_{\mathcal
E})\mapsto\ker(\dbar_{\mathcal E})$ are mutually inverse
equivalences of suitably defined categories.
\end{remark}

By the very definition, the kernel of the canonical $\dbar$-operator
on ${\mathcal T}_X^{\mathbb C}/{\mathcal P}$ coincides with
$\left({\mathcal T}_X^{\mathbb C}/{\mathcal P}\right)^{\mathcal P}$
(the subsheaf of ${\mathcal P}$-invariant sections, see Example
\ref{hol vector fields}). We denote this subsheaf by ${\mathcal
T}_{X/{\mathcal P}}$.

The duality pairing restricts to a non-degenerate ${\mathcal
O}_{X/{\mathcal P}}$-bilinear pairing between $\Omega^1_{X/{\mathcal
P}}$ and ${\mathcal T}_{X/{\mathcal P}}$ giving rise to a faithful
action of ${\mathcal T}_{X/{\mathcal P}}$ on ${\mathcal
O}_{X/{\mathcal P}}$ by derivations by the usual formula
$\xi(f)=\iota_\xi df$, for $\xi\in{\mathcal T}_{X/{\mathcal P}}$ and
$f\in{\mathcal O}_{X/{\mathcal P}}$.

Let $\pr_i:X\times X\to X$ denote the projection on the
$i^{\text{th}}$ factor and let $\Delta_X : X\to X\times X$ denote
the diagonal embeding. The latter satisfies $\Delta^*({\mathcal
P}\times{\mathcal P})={\mathcal P}$. Therefore, the induced map
$\Delta_X^*: {\mathcal O}_{X\times X}\to{\mathcal O}_X$ satisfies
\[
\im(\left.\Delta_X^*\right|_{{\mathcal O}_{{X\times X}/{\mathcal
P}\times{\mathcal P}}})\subset{\mathcal O}_{X/{\mathcal P}} \ .
\]
Let $\Delta_{X/{\mathcal P}}^*:=\left.\Delta_X^*\right|_{{\mathcal
O}_{{X\times X}/{\mathcal P}\times{\mathcal P}}}$, ${\mathcal
I}_{X/{\mathcal P}}:=\ker(\Delta_{X/{\mathcal P}}^*)$.

For a locally-free ${\mathcal O}_{X/{\mathcal P}}$-module of finite
rank ${\mathcal E}$ let
\[
{\mathcal J}^k({\mathcal E}):=(\pr_1)_*\left({\mathcal O}_{X\times
X/{\mathcal P}\times{\mathcal P}}/{\mathcal I}_{X/{\mathcal
P}}^{k+1}\otimes_{\pr_2^{-1}{\mathcal O}_{X/{\mathcal
P}}}\pr_2^{-1}{\mathcal E}\right) \ ,
\]
let ${\mathcal J}^k_{X/{\mathcal P}}:= {\mathcal J}^k({\mathcal
O}_{X/{\mathcal P}})$. It is clear from the above definition that
${\mathcal J}^k_{X/{\mathcal P}}$ is, in a natural way, a
commutative algebra and ${\mathcal J}^k({\mathcal E})$ is a
${\mathcal J}^k_{X/{\mathcal P}}$-module.

We regard ${\mathcal J}^k({\mathcal E})$ as ${\mathcal
O}_{X/{\mathcal P}}$-modules via the pull-back map
$\pr_1^*:{\mathcal O}_{X/{\mathcal P}}\to(\pr_1)_*{\mathcal
O}_{{X\times X}/{\mathcal P}\times{\mathcal P}}$ (the restriction to
${\mathcal O}_{X/{\mathcal P}}$ of the map $\pr_1^*:{\mathcal
O}_X\to(\pr_1)_*{\mathcal O}_{X\times X}$) with the quotient map
$(\pr_1)_*{\mathcal O}_{{X\times X}/{\mathcal P}\times{\mathcal
P}}\to{\mathcal J}_{X/{\mathcal P}}^k$.

For $0\leq k\leq l$ the inclusion ${\mathcal I}_{X/{\mathcal
P}}^{l+1}\to{\mathcal I}_{X/{\mathcal P}}^{k+1}$ induces the
surjective map $\pi_{l,k}:{\mathcal J}^l({\mathcal E})\to{\mathcal
J}^k({\mathcal E})$. The sheaves ${\mathcal J}^k({\mathcal E})$,
$k=0,1,\ldots$ together with the maps $\pi_{l,k}$, $k\leq l$ form an
inverse system. Let ${\mathcal J}({\mathcal E})={\mathcal
J}^\infty({\mathcal E}):=\underset{\longleftarrow}{\lim}{\mathcal
J}^k({\mathcal E})$. Thus, ${\mathcal J}({\mathcal E})$ carries a
natural topology.

Let $j^k: {\mathcal E}\to{\mathcal J}^k({\mathcal E})$ denote the
map $e\mapsto 1\otimes e$,
$j^\infty:=\underset{\longleftarrow}{\lim}j^k$.

Let
\begin{multline*}
d_1 : {\mathcal O}_{{X\times X}/{\mathcal P}\times{\mathcal
P}}\otimes_{\pr_2^{-1}{\mathcal O}_{X/{\mathcal
P}}}\pr_2^{-1}{\mathcal E} \to \\
\to \pr_1^{-1}\Omega^1_{X/{\mathcal P}}\otimes_{\pr_1^{-1}{\mathcal
O}_{X/{\mathcal P}}}{\mathcal O}_{{X\times X}/{\mathcal
P}\times{\mathcal P}}\otimes_{\pr_2^{-1}{\mathcal O}_{X/{\mathcal
P}}}\pr_2^{-1}{\mathcal E}
\end{multline*}
denote the exterior derivative along the first factor. It satisfies
\begin{multline*}
d_1({\mathcal I}_{X/{\mathcal P}}^{k+1}\otimes_{\pr_2^{-1}{\mathcal
O}_{X/{\mathcal P}}}\pr_2^{-1}{\mathcal E})\subset
\\
\pr_1^{-1}\Omega^1_X\otimes_{\pr_1^{-1}{\mathcal O}_{X/{\mathcal
P}}}{\mathcal I}_{X/{\mathcal P}}^k\otimes_{\pr_2^{-1}{\mathcal
O}_{X/{\mathcal P}}}\pr_2^{-1}{\mathcal E}
\end{multline*}
for each $k$ and, therefore, induces the map
\[
d_1^{(k)} : {\mathcal J}^k({\mathcal E})\to\Omega^1_{X/{\mathcal
P}}\otimes_{{\mathcal O}_{X/{\mathcal P}}}{\mathcal
J}^{k-1}({\mathcal E})
\]
The maps $d_1^{(k)}$ for different values of $k$ are compatible with
the maps $\pi_{l,k}$ giving rise to the \emph{canonical flat
connection}
\[
\nabla^{can}_{\mathcal E} : {\mathcal J}({\mathcal
E})\to\Omega^1_{X/{\mathcal P}}\otimes_{{\mathcal O}_{X/{\mathcal
P}}}{\mathcal J}({\mathcal E})
\]
which extends to the flat connection
\[
\nabla^{can}_{\mathcal E} : {\mathcal O}_X\otimes_{{\mathcal
O}_{X/{\mathcal P}}}{\mathcal J}({\mathcal
E})\to\Omega^1_X\otimes_{{\mathcal O}_X}({\mathcal
O}_X\otimes_{{\mathcal O}_{X/{\mathcal P}}}{\mathcal J}({\mathcal
E}))
\]
Here and below by abuse of notation we write
$(\bullet)\otimes_{{\mathcal O}_{X/{\mathcal P}}}{\mathcal
J}_{X/{\mathcal P}}({\mathcal E})$ for
$\underset{\longleftarrow}{\lim}(\bullet)\otimes_{{\mathcal
O}_{X/{\mathcal P}}}{\mathcal J}^k({\mathcal E})$.

\subsection{Hochschild cochains}
Suppose that ${\mathcal A}$ is an ${\mathcal O}_{X/{\mathcal
P}}$-Azumaya algebra , i.e. a sheaf of algebras on $X$ locally
isomorphic to $\Mat_n({\mathcal O}_{X/{\mathcal P}})$.

For $n>0$ let $C^n({\mathcal A})$ denote the sheaf of
multidifferential operators ${\mathcal A}^{\times n}\to{\mathcal
A}$; Let $C^0({\mathcal A})={\mathcal A}$. The subsheaf of
normalized cochains $\overline C^n({\mathcal A})$ consists of those
multidifferential operators which yield zero whenever one of the
arguments is in ${\mathbb C}\cdot 1\subset{\mathcal A}$.

With the Gerstenhaber bracket $[\ ,\ ]$ and the Hochschild
differential, denoted $\delta$, defined in the standard fashion,
$C^\bullet({\mathcal A})[1]$ and $\overline C^\bullet({\mathcal
A})[1]$ are DGLA and the inclusion $\overline C^\bullet({\mathcal
A})[1]\to C^\bullet({\mathcal A})[1]$ is a quasi-isomorphism of
such.

For $n>0$ let $C^n({\mathcal O}_X\otimes_{{\mathcal O}_{X/{\mathcal
P}}}{\mathcal J}({\mathcal A}))$ denote the sheaf of continuous
${\mathcal O}_X$-multilinear maps $({\mathcal O}_X\otimes_{{\mathcal
O}_{X/{\mathcal P}}}{\mathcal J}({\mathcal A}))^{\times
n}\to{\mathcal O}_X\otimes_{{\mathcal O}_{X/{\mathcal P}}}{\mathcal
J}({\mathcal A})$. The canonical flat connection
$\nabla^{can}_{\mathcal A}$ on ${\mathcal O}_X\otimes_{{\mathcal
O}_{X/{\mathcal P}}}{\mathcal J}({\mathcal A})$ induces the flat
connection, still denoted $\nabla^{can}_{\mathcal A}$, on
$C^n({\mathcal O}_X\otimes_{{\mathcal O}_{X/{\mathcal P}}}{\mathcal
J}({\mathcal A}))$. Equipped with the Gerstenhaber bracket and the
Hochschild differential $\delta$, $C^\bullet({\mathcal
O}_X\otimes_{{\mathcal O}_{X/{\mathcal P}}}{\mathcal J}({\mathcal
A}))[1]$ is a DGLA. Just as in the case ${\mathcal P}=0$ (see
\ref{ssec:cochains}) we have the DGLA
$\Omega^\bullet_X\otimes_{{\mathcal O}_X}C^\bullet({\mathcal
O}_X\otimes_{{\mathcal O}_{X/{\mathcal P}}}{\mathcal J}({\mathcal
A}))[1]$ with the differential $\nabla^{can}+\delta$.

\begin{lemma}
The de Rham complex $\DR(C^n({\mathcal O}_X\otimes_{{\mathcal
O}_{X/{\mathcal P}}}{\mathcal J}({\mathcal
A}))):=(\Omega^\bullet_X\otimes_{{\mathcal O}_X}C^n({\mathcal
O}_X\otimes_{{\mathcal O}_{X/{\mathcal P}}}{\mathcal J}({\mathcal
A})),\nabla^{can}_{\mathcal A})$ satisfies
\begin{enumerate}
\item $H^i\DR(C^n({\mathcal O}_X\otimes_{{\mathcal
O}_{X/{\mathcal P}}}{\mathcal J}({\mathcal A})))=0$ for $i\neq 0$
\item The map $j^\infty : C^n({\mathcal A})\to C^n({\mathcal O}_X\otimes_{{\mathcal O}_{X/{\mathcal P}}}{\mathcal J}({\mathcal
A}))$ is an isomorphism onto $H^0\DR(C^n({\mathcal
O}_X\otimes_{{\mathcal O}_{X/{\mathcal P}}}{\mathcal J}({\mathcal
A})))$.
\end{enumerate}
\end{lemma}
\begin{cor}
The map $j^\infty:C^\bullet({\mathcal
A})\to\Omega^\bullet_X\otimes_{{\mathcal O}_X} C^\bullet({\mathcal
O}_X\otimes_{{\mathcal O}_{X/{\mathcal P}}}{\mathcal J}({\mathcal
A}))[1]$ is a quasi-isomorphism of DGLA.
\end{cor}

\subsection{Azumaya} Suppose that ${\mathcal A}$ is an ${\mathcal
O}_{X/{\mathcal P}}$-Azumaya algebra.

The sheaves ${\mathcal J}({\mathcal A})$ and ${\mathcal
A}\otimes_{{\mathcal O}_{X/{\mathcal P}}}{\mathcal J}_{X/{\mathcal
P}}$ have canonical structures of sheaves of (central) ${\mathcal
J}_{X/{\mathcal P}}$-algebras locally isomorphic to
$\Mat_n({\mathcal J}_{X/{\mathcal P}})$ and come equipped with
projections to ${\mathcal A}$.

Let $\shIsom_0({\mathcal O}_X\otimes_{{\mathcal O}_{X/{\mathcal P}}}
{\mathcal A}\otimes_{{\mathcal O}_{X/{\mathcal P}}}{\mathcal
J}_{X/{\mathcal P}},{\mathcal O}_X\otimes_{{\mathcal O}_{X/{\mathcal
P}}}{\mathcal J}({\mathcal A}))$ denote the sheaf of (locally
defined) ${\mathcal O}_X\otimes_{{\mathcal O}_{X/{\mathcal
P}}}{\mathcal J}_{X/{\mathcal P}}$-algebra isomorphisms ${\mathcal
O}_X\otimes_{{\mathcal O}_{X/{\mathcal P}}}{\mathcal
A}\otimes_{{\mathcal O}_{X/{\mathcal P}}}{\mathcal J}_{X/{\mathcal
P}}\to{\mathcal O}_X\otimes_{{\mathcal O}_{X/{\mathcal P}}}{\mathcal
J}({\mathcal A})$ which induce the identity map on ${\mathcal
O}_X\otimes_{{\mathcal O}_{X/{\mathcal P}}}{\mathcal A}$. Let
$\shAut_0({\mathcal O}_X\otimes_{{\mathcal O}_{X/{\mathcal
P}}}{\mathcal A}\otimes_{{\mathcal O}_{X/{\mathcal P}}}{\mathcal
J}_{X/{\mathcal P}})$ denote the sheaf of (locally defined)
${\mathcal O}_X\otimes_{{\mathcal O}_{X/{\mathcal P}}}{\mathcal
J}_{X/{\mathcal P}}$-algebra automorphisms of ${\mathcal
O}_X\otimes_{{\mathcal O}_{X/{\mathcal P}}}{\mathcal
A}\otimes_{{\mathcal O}_{X/{\mathcal P}}}{\mathcal J}_{X/{\mathcal
P}}$ which induce the identity map on ${\mathcal
O}_X\otimes_{{\mathcal O}_{X/{\mathcal P}}}{\mathcal A}$.

Just as in Corollary \ref{trivial torsor} we may conclude that
\begin{multline*}
\Isom_0({\mathcal O}_X\otimes_{{\mathcal O}_{X/{\mathcal P}}}
{\mathcal A}\otimes_{{\mathcal O}_{X/{\mathcal P}}}{\mathcal
J}_{X/{\mathcal P}},{\mathcal O}_X\otimes_{{\mathcal O}_{X/{\mathcal
P}}}{\mathcal J}({\mathcal A})):= \\
\Gamma(X;\shIsom_0({\mathcal
O}_X\otimes_{{\mathcal O}_{X/{\mathcal P}}} {\mathcal
A}\otimes_{{\mathcal O}_{X/{\mathcal P}}}{\mathcal J}_{X/{\mathcal
P}},{\mathcal O}_X\otimes_{{\mathcal O}_{X/{\mathcal P}}}{\mathcal
J}({\mathcal A})))
\end{multline*}
is non-empty.

A choice of $\sigma\in\Isom_0({\mathcal O}_X\otimes_{{\mathcal
O}_{X/{\mathcal P}}} {\mathcal A}\otimes_{{\mathcal O}_{X/{\mathcal
P}}}{\mathcal J}_{X/{\mathcal P}},{\mathcal O}_X\otimes_{{\mathcal
O}_{X/{\mathcal P}}}{\mathcal J}({\mathcal A}))$ and
$\nabla\in{\mathcal C}({\mathcal O}_X\otimes_{{\mathcal
O}_{X/{\mathcal P}}}{\mathcal A})$ give rise to a unique
$F\in\Gamma(X;\Omega^1_X\otimes_{{\mathcal O}_{X/{\mathcal
P}}}{\mathcal A}_0\otimes_{{\mathcal O}_{X/{\mathcal P}}}{\mathcal
J}_{X/{\mathcal P}})$ and $\theta
\in\Gamma(X;\Omega^2_X\otimes_{{\mathcal O}_{X/{\mathcal
P}}}{\mathcal A}_0)$ such that $\nabla^2 =\ad \theta$ and the
equation \eqref{3.0.4} holds. Such a $\sigma$ provides us with the
isomorphism of DGLA
\begin{multline}\label{quism of DGLA with pol}
\sigma_*: \Omega^\bullet_X\otimes_{{\mathcal
O}_X}C^\bullet({\mathcal O}_X\otimes_{{\mathcal O}_{X/{\mathcal
P}}}{\mathcal A}\otimes_{{\mathcal O}_{X/{\mathcal P}}}{\mathcal
J}_{X/{\mathcal
P}})[1]\to \\
\to\Omega^\bullet_X\otimes_{{\mathcal O}_X} C^\bullet({\mathcal
O}_X\otimes_{{\mathcal O}_{X/{\mathcal P}}}{\mathcal J}({\mathcal
A}))[1]
\end{multline}
where the latter is equipped with the differential
$\nabla^{can}_{\mathcal A}+\delta$ and the former is equipped with
the differential $\nabla\otimes\id+\id\otimes\nabla^{can}_{\mathcal
O}+\ad(F)+\delta$.

The operator $\exp(\iota_F)$ is an automorphism of the graded Lie
algebra $\Omega^\bullet_X\otimes_{{\mathcal O}_X}C^\bullet({\mathcal
O}_X\otimes_{{\mathcal O}_{X/{\mathcal P}}}{\mathcal
A}\otimes_{{\mathcal O}_{X/{\mathcal P}}}{\mathcal J}_{X/{\mathcal
P}})[1]$. It does not commute with the differential
$\nabla\otimes\id + \id\otimes\nabla^{can}_{\mathcal O} + \ad(F) +
\delta$. Instead, the formula of Lemma \ref{Adiota} holds. Hence,
the composition $\sigma_*\circ\exp(\iota_F)$ is a quasi-isomorphism
of DGLA as in \eqref{quism of DGLA with pol} but with the source
equipped with the differential $\nabla\otimes\id +
\id\otimes\nabla^{can}_{\mathcal O} + \delta + \iota_\theta$.

The cotrace map
\begin{multline*}
\cotr : \overline C^\bullet({\mathcal O}_X\otimes_{{\mathcal
O}_{X/{\mathcal P}}}{\mathcal J}_{X/{\mathcal P}})[1]\to \\
\to C^\bullet({\mathcal O}_X\otimes_{{\mathcal O}_{X/{\mathcal
P}}}{\mathcal A}\otimes_{{\mathcal O}_{X/{\mathcal P}}}{\mathcal
J}_{X/{\mathcal P}})[1]
\end{multline*}
defined as in \eqref{define cotrace} gives rise to the
quasi-isomorphism of DGLA
\begin{multline*}
\id\otimes\cotr : \Omega^\bullet_X\otimes_{{\mathcal O}_X}\overline
C^\bullet({\mathcal O}_X\otimes_{{\mathcal O}_{X/{\mathcal
P}}}{\mathcal J}_{X/{\mathcal P}})[1]\to \\
\to \Omega^\bullet_X\otimes_{{\mathcal O}_X} C^\bullet({\mathcal
O}_X\otimes_{{\mathcal O}_{X/{\mathcal P}}}{\mathcal
A}\otimes_{{\mathcal O}_{X/{\mathcal P}}}{\mathcal J}_{X/{\mathcal
P}})[1]
\end{multline*}
where the source (respectively, the target) is equipped with the
differential $\nabla^{can}_{\mathcal O} + \delta$ (respectively,
$\nabla\otimes\id+\id\otimes\nabla^{can}_{\mathcal O} + \delta +
\iota_\theta$).

The proof of \ref{prop:independence of choices} shows that the image
of the composition
$\sigma_*\circ\exp(\iota_{F})\circ(\id\otimes\cotr)$ in the derived
category does not depend on the choices made. We summarize the above
in the following theorem.

\begin{thm}\label{quism jets pol}
Suppose that ${\mathcal P}$ is an integrable (complex) distribution
on $X$ and ${\mathcal A}$ is an ${\mathcal O}_{X/{\mathcal
P}}$-Azumaya algebra. There is a canonical isomorphism in the
derived category of DGLA $\Omega^\bullet_X\otimes_{{\mathcal O}_X}
\overline C^\bullet({\mathcal O}_X\otimes_{{\mathcal O}_{X/{\mathcal
P}}}{\mathcal J}_{X/{\mathcal
P}})[1]\cong\Omega^\bullet_X\otimes_{{\mathcal
O}_X}C^\bullet({\mathcal O}_X\otimes_{{\mathcal O}_{X/{\mathcal P}}}
{\mathcal J}({\mathcal A}))[1]$.
\end{thm}
\begin{cor}
Under the assumptions of Theorem \ref{quism jets pol}, there is a
canonical isomorphism in the derived category of DGLA $\overline
C^\bullet({\mathcal O}_{X/{\mathcal P}})\cong C^\bullet({\mathcal
A})$.
\end{cor}

\end{document}